\documentclass[11pt]{article}
\usepackage{graphicx}
\usepackage{color}
\usepackage{amsmath}
\usepackage{amssymb}
\usepackage{amscd}
\usepackage{bbm}

\newcommand{\R}{\mathbb{R}}
\newcommand{\inr}[1]{\bigl< #1 \bigr>}

\newcommand{\sgn}{{\rm sgn}}

\newcommand{\E}{\mathbb{E}}

\newcommand{\eps}{\varepsilon}

\newcommand{\cF}{{\cal F}}

\newcommand{\cG}{{\cal G}}
\newcommand{\cX}{{\cal X}}
\newcommand{\cN}{{\cal N}}

\newcommand{\vertiii}[1]{{\left\vert\kern-0.25ex\left\vert\kern-0.25ex\left\vert #1\right\vert\kern-0.25ex\right\vert\kern-0.25ex\right\vert}}

\newtheorem{Theorem}{Theorem}[section]
\newtheorem{Lemma}[Theorem]{Lemma}

\newtheorem{Definition}[Theorem]{Definition}

\newtheorem{Corollary}[Theorem]{Corollary}
\newtheorem{Remark}[Theorem]{Remark}

\numberwithin{equation}{section}

\newcommand{\norm}[1]{\left\|#1\right\|}%

\DeclareMathOperator*{\argmin}{argmin}

\def \proof {\noindent {\bf Proof.}\ \ }

\def \endproof
{{\mbox{}\nolinebreak\hfill\rule{2mm}{2mm}\par\medbreak}}
\def\IND{\mathbbm{1}}

\begin{document}
\title{{Sparse recovery under weak moment assumptions}}
\author{Guillaume Lecu\'e${}^{1,3}$  \and Shahar Mendelson${}^{2,4,5}$}

\footnotetext[1]{CNRS, CMAP, Ecole Polytechnique, 91120 Palaiseau, France.}
\footnotetext[2]{Department of Mathematics, Technion, I.I.T, Haifa
32000, Israel.}
 \footnotetext[3] {Email:
guillaume.lecue@cmap.polytechnique.fr }
\footnotetext[4] {Email:
shahar@tx.technion.ac.il}
\footnotetext[5]{Supported by the Mathematical Sciences Institute -- The Australian National University and by ISF grant 900/10.}

\maketitle

\begin{abstract}
  We prove that iid random vectors that satisfy a rather weak moment
  assumption can be used as measurement vectors in Compressed Sensing, and the number of measurements required for exact reconstruction is the same as the best possible estimate -- exhibited by a random Gaussian matrix. We then show that this moment condition is necessary, up to a $\log \log $ factor. In addition, we explore the Compatibility Condition and the Restricted Eigenvalue Condition in the noisy setup, as well as properties of neighbourly random polytopes.
\end{abstract}

\section{Introduction and main results}
\label{sec:intro}
Data acquisition is an important task in diverse fields such as mobile communications, medical imaging, radar detection and others, making the design of efficient data acquisition processes a problem of obvious significance.

The core issue in data acquisition is retaining all the valuable information at one's disposal, while keeping the `acquisition cost' as low as possible. And while there are several ways of defining that cost, depending on the problem (storage, time, financial cost, etc.), the common denominator of being `cost effective' is ensuring the quality of the data while keeping the number of measurements as small as possible.

The rapidly growing area of {\it Compressed Sensing} studies `economical' data acquisition
processes. We refer the reader to \cite{MR2230846,MR2241189} and to
the book \cite{MR3100033} for more information on the origins of Compressed Sensing and a survey of the progress that has been made in the area in recent years.

At the heart of Compressed Sensing is a simple idea that has been a recurring theme in Mathematics and Statistics: while complex objects (in this case, data), live in high-dimensional spaces, they can be described effectively using low-dimensional, approximating structures; moreover, randomness may be used to expose these low-dimensional structures. Of course, unlike more theoretical applications of this idea, identifying the low-dimensional structures in the context of Compressed Sensing must be robust and efficient, otherwise, such procedures will be of little practical use.

In the standard Compressed Sensing setup, one observes linear measurements
$y_i=\inr{X_i,x_0}$, $i=1,...,N$ of an unknown vector $x_0 \in \R^n$. To make the data acquisition process `cost-effective', the number of measurements $N$ is assumed to be much smaller than the dimension $n$, and the goal is to identify $x_0$ using those measurements.

Because the resulting system of equations is under-determined, there is no hope, in general, of identifying $x_0$. However, if
$x_0$ is believed to be well approximated by a low-dimensional structure, for example, if $x_0$ is supported on at most $s$ coordinates for some $s \leq N$, the problem becomes more feasible.

Let $(f_1,...,f_N)$ be the canonical basis of $\R^N$ (we will later
use $(e_1,\ldots,e_n)$ to denote the canonical basis of $\R^n$) and consider the matrix
$$
\Gamma=\frac{1}{\sqrt{N}}\sum_{i=1}^N  \inr{X_i,\cdot}f_i,
$$
called the {\it measurement matrix}. One possible recovery procedure is $\ell_0$-minimization, in which one selects a vector $t\in\R^n$ that has the shortest support among all vectors satisfying $\Gamma t=\Gamma x_0$. Unfortunately, $\ell_0$ minimization is known to be NP-hard in general (see \cite{MR1320206} or Theorem~2.17 in \cite{MR3100033}). Thus, even without analyzing if and when $\ell_0$-minimization actually recovers $x_0$, it is obvious that a more computationally reasonable procedure has to be found.

Fortunately, efficient  procedures have been used since the seventies in
geophysics (see, for instance \cite{claerbout},\cite{tbm}, 
\cite{MR857796}, and Logan's
Ph.D. thesis \cite{logan}). Those procedures are based on
$\ell_1$-minimization for which  early theoretical works can be found in \cite{MR1154788} and \cite{MR1854649}.

In particular, \textit{Basis Pursuit} is  a convex relaxation of $\ell_0$-minimization, and since it can be recast as a linear program (see, e.g., Chapter~15 in \cite{MR3100033}), it is far more reasonable than $\ell_0$-minimization from the computational viewpoint.
\begin{Definition} \label{def:basis-pursuit}
Given the measurement matrix $\Gamma$ and the measurements $\Gamma x_0=(\inr{X_i,x_0})_{i=1}^N$, Basis Pursuit returns a vector $\hat{x}$ that satisfies
\begin{equation}\label{eq:basis-pursuit}
  \hat x\in{\rm argmin}\big(\|t\|_1:\Gamma t=\Gamma x_0\big).
\end{equation}
\end{Definition}

Since one may solve this minimization problem effectively, the focus may be shifted to the quality of the solution: whether one can identify measurement vectors $X_1,....,X_N$ for which \eqref{eq:basis-pursuit} has a unique solution, which is $x_0$ itself, for any $x_0$ that is $s$-sparse (i.e. supported on at most $s$ coordinates).
\begin{Definition}\label{def:ER}
Let $\Sigma_s$ be the set of all $s-$sparse vectors in $\R^n$. An ${N\times n}$ matrix $\Gamma$ satisfies the \textbf{exact reconstruction property of order $s$} if for every $x_0\in\Sigma_s$,
  \begin{equation}
    \argmin\big(\norm{t}_1 : \Gamma t=\Gamma x_0 \big)=\{x_0\}. \tag{ER(s)}
  \end{equation}
\end{Definition}

It follows from Proposition~2.2.18 in \cite{MR3113826} that if
$\Gamma$ satisfies ER($s$) then necessarily the number of measurements (rows)
is at least $N\geq c_0 s \log\big(en/s\big)$, where $c_0$ is a suitable absolute constant.
On the other hand, there are constructions of (random) matrices $\Gamma$ that satisfy ER($s$) with $N$ proportional to $s\log\big(en/s\big)$. From here on and with a minor abuse of notation, we will refer to $s\log(en/s)$ as the {\it optimal number of measurements} and ignore the exact dependence on the constant $c_0$.

Unfortunately, the only matrices that are known to satisfy the reconstruction property with an optimal number of measurements are random -- which is not surprising, as randomness is one of the most effective tools in exposing low-dimensional, approximating structures.  A typical example of an `optimal matrix' is the Gaussian matrix, which has independent standard normal random variables as
entries. Other examples of optimal measurement matrices are $\Gamma=N^{-1/2}\sum_{i=1}^N \inr{X_i,\cdot}f_i$ where $X_1,...,X_N$ are independent, isotropic and $L$-subgaussian random
vectors:
\begin{Definition} \label{def:subgauss}
A symmetric random vector $X \in \R^n$ is isotropic if for every $t \in \R^n$, $\E \inr{X,t}^2 = \|t\|_2^2$; it is $L$-subgaussian if for every $t \in \R^n$ and every $p \geq 2$, $\|\inr{X,t}\|_{L_p} \leq L \sqrt{p} \|\inr{X,t}\|_{L_2}$.
\end{Definition}

The optimal behaviour of isotropic, $L$-subgaussian matrix ensembles and other ensembles like it, occurs because a typical matrix acts on $\Sigma_s$ in an isomorphic way when $N\geq c_1s \log(en/s)$, and in the $L$-subgaussian case, $c_1$ is a constant that depends only on $L$. In Compressed Sensing literature, this isomorphic behaviour is called the {\it Restricted Isometry property} (RIP) (see, for example \cite{MR2236170,MR2300700,MR2453368}): A matrix $\Gamma$ satisfies the RIP in $\Sigma_s$ with constant $0<\delta<1$, if for every $t \in \Sigma_s$,
\begin{equation}\label{eq:RIP}
(1-\delta) \|t\|_2 \leq \|\Gamma t\|_2 \leq (1+\delta)\|t\|_2.
\end{equation}
It is straightforward to show that if $\Gamma$ satisfies the RIP in $\Sigma_{2s}$ for a sufficiently small constant $\delta$, then it has the exact reconstruction property of order $s$ (see, e.g. \cite{MR2230846,MR2300700,C:08}).

The standard proof of the RIP for subgaussian ensembles is based on the rapid tail decay of linear functionals $\inr{X,t}$. Thus, it seemed natural to ask whether the RIP holds even when linear functionals exhibit a slower decay -- for example, when $X$ is
$L$-subexponential -- that is, when linear functionals only satisfy that $\|\inr{X,t}\|_{L_p} \leq L p \|\inr{X,t}\|_{L_2}$ for every $t \in \R^n$ and every $p \geq 2$.

Proving the RIP for subexponential ensembles is a much harder task than for subgaussian ensembles (see, e.g. \cite{MR2796091}). Moreover, the RIP does not exhibit the same optimal quantitative behaviour as in the Gaussian case: it holds with high probability only when $N\geq  c_2(L)s \log^{2}(en/s)$, and this estimate cannot be improved, as can be seen when $X$ has independent, symmetric exponential random variables as coordinates \cite{MR2796091}.

Although the RIP need not be true for isotropic $L$-subexponential ensemble using the optimal number of measurements, results in \cite{MR2829871} (see Theorem~7.3 there) and in \cite{MR3134335} show that exact reconstruction can still be achieved by such an ensemble and {\it with the optimal number of measurements}. This opens the door to an intriguing question: whether considerably weaker assumptions on the measurement vector may still lead to Exact Reconstruction even when the RIP fails.

The main result presented here does just that, using the {\it small-ball method} introduced in \cite{Shahar-COLT,Shahar-Gelfand}.

\begin{Definition}\label{def:small_ball_d_sparse}
 A random vector $X$ satisfies the \textbf{small-ball condition} in the set $\Sigma_s$ with constants $u,\beta>0$ if for every $t \in \Sigma_s$,
 $$
 P\big(|\inr{X,t}|>u \norm{t}_2\big)\geq\beta.
 $$
\end{Definition}
The small-ball condition is a rather minimal assumption on the measurement vector and is satisfied in fairly general situations for values of $u$ and $\beta$ that are suitable constants, independent of the dimension $n$.

Under some normalization (like isotropicity), a small-ball condition is an immediate outcome of the Paley-Zygmund inequality (see, e.g. \cite{MR1666908}) and moment equivalence. For example, in the following cases a small-ball condition holds with constants that depend only on $\kappa_0$ (and on $\eps$ for the first case); the straightforward proof may be found in  \cite{Shahar-COLT}.
\begin{description}
\item{$\bullet$} $X$ is isotropic and for every  $t\in\Sigma_s$, $\norm{\inr{X,t}}_{L_{2+\eps}}\leq \kappa_0 \norm{\inr{X,t}}_{L_2}$ for some $\eps>0$;
\item{$\bullet$}  $X$ is isotropic and for every  $t\in\Sigma_s$,
  $\norm{\inr{X,t}}_{L_2}\leq \kappa_0 \norm{\inr{X,t}}_{L_1}$.
\end{description}

Because the small-ball condition means that marginals of $X$ do not assign too much weight close to $0$, it may hold even without integrability (and in particular, $X$ need not have a covariance matrix). One such example is a random vector whose coordinates are independent random variables that are
absolutely continuous with respect to the Lebesgue measure and with a
density  almost surely bounded by $\kappa_0$. Indeed, as noted in Theorem~1.2 from \cite{RV13}, for every $t \in \R^n$, $\inr{X,t/\|t\|_2}$
has a density that is almost surely bounded by $\sqrt{2}\kappa_0$.  In
particular, $P\big(|\inr{X,t}|\geq (4\sqrt{2}\kappa_0)^{-1}\|t\|_2\big)\geq
1/2$ and $X$ satisfies the small ball condition with
$u=(4\sqrt{2}\kappa_0)^{-1}$ and $\beta=1/2$. The estimate on the density of $\inr{X,t/\|t\|_2}$ follows by combining a result due to B. Rogozin \cite{MR890930} on the maximal value of a
convolution product of densities, and a result due to K. Ball
\cite{MR840631}, on the maximal volume of a section of the cube
$[-1/2,1/2]^n$.

\vskip0.5cm

Our first result shows that a combination of the small-ball condition and a weak moment assumption suffices to ensure the exact reconstruction property with the optimal number of measurements.

\vspace{0.4cm}

\noindent {\bf Theorem A.} \textit{There exist absolute constants
  $c_0$, $c_1$ and $c_2$ and for every $\alpha \geq 1/2$ there exists
  a constant $c_3(\alpha)$ that depends only on $\alpha$ for which the
  following holds. Let $X=(x_i)_{i=1}^n$ be a random vector on $\R^n$
  (with potentially dependent coordinates). Assume that
\begin{enumerate}
\item[1.] there are $\kappa_1,\kappa_2,w>1$ that satisfy that for every $1 \leq j \leq n$,
$\|x_j\|_{L_2} =1$ and, for every $4\leq p\leq 2\kappa_2\log(wn)$, $\|x_j\|_{L_p} \leq \kappa_1 p^\alpha$.
\item[2.] $X$ satisfies the small ball condition in $\Sigma_s$ with constants $u$ and $\beta$.
\end{enumerate}
If
\begin{equation*}
N \geq c_0\max\left\{s \log \Big(\frac{en}{s}\Big), (c_3(\alpha)\kappa_1^2)^2 (\kappa_2\log(wn))^{\max\{4\alpha-1,1\}} \right\},
\end{equation*}
and $X_1,...,X_N$ are independent copies of $X$, then,  with probability at least
$$
1-2\exp(-c_1\beta^2N)-1/w^{\kappa_2} n^{\kappa_2-1},
$$
$\Gamma=N^{-1/2}\sum_{i=1}^N \inr{X_i,\cdot}f_i$ satisfies the exact reconstruction property in $\Sigma_{s_1}$ for $s_1=c_2u^2 \beta s$.
}
\vskip0.5cm
An immediate outcome of Theorem~A is the following:
\begin{description}
\item{$\bullet$} Let $x$ be a centered random variable that has
  variance $1$ and for which $\|x\|_{L_p} \leq c\sqrt{p}$ for $1\leq
  p\leq 2\log n$. If $X$ has independent coordinates distributed as
  $x$, then the corresponding matrix $\Gamma$ with $N\geq c_1s
  \log(en/s)$ rows can be used as a measurement matrix and recover any
  $s$-sparse vector with large probability.
\end{description}
It is relatively straightforward to derive many other results of a
similar flavour, leading to random ensembles that satisfy the exact
reconstruction property with the optimal number of measurements.

\begin{Remark}
Our focus is on measurement matrices with independent rows, that satisfy conditions of a stochastic nature -- they have i.i.d. rows.  Other
types of measurement matrices that have some \textit{structure} have also been
used in Compressed Sensing. One notable example is a random Fourier measurement matrix, obtained by randomly selecting rows from
the discrete Fourier matrix (see, e.g.  \cite{MR2300700}, \cite{MR2417886} or Chapter~12 in
\cite{MR3100033}). 

One may wonder if the small-ball condition is satisfied for more structured matrices, as the argument we use here does not extend immediately to such cases.
And, indeed, for structured ensembles one may encounter a different situation: a small-ball condition that is not uniform, in the sense that
the constants $u$ and $\beta$ from Definition~\ref{def:small_ball_d_sparse} are direction-dependent. Moreover, in some cases, the known estimates on these constants are far from what is expected.

Results of the same flavour of Theorem~A may follow from a `good enough' small-ball condition, even if it is not uniform, by slightly modifying the argument we use here. However,
obtaining a satisfactory `non-uniform' small-ball condition is a different story. For example, in the Fourier case, such an estimate is likely to require quantitative extensions of the Littlewood-Paley theory -- a worthy challenge in its own right, and one which goes far beyond the goals of this article.
\end{Remark}

\vskip0.3cm

Just as noted for subexponential ensembles, Theorem~A cannot be proved using an RIP-based argument. A key ingredient in the proof is the following observation:

\vspace{0.5cm}
\noindent {\bf Theorem B.} \textit{Let $\Gamma:\R^n\mapsto\R^N$ and
  denote by  $(e_1,\ldots,e_n)$ the canonical basis of $\R^n$. Assume
  that:
\begin{enumerate}
\item[a)] for every $x \in \Sigma_s$, $\norm{\Gamma x}_2\geq c_0 \norm{x}_2$, and
\item[b)] for every $j\in\{1,\ldots,n\}$, $\norm{\Gamma e_j}_2\leq c_1$.
\end{enumerate}
Setting $s_1=\big\lfloor(c_0^2(s-1))/(4c_1^2)\big\rfloor-1$, $\Gamma$ satisfies the exact reconstruction property in $\Sigma_{s_1}$.  }

\vspace{0.5cm}

Compared with the RIP, conditions a) and b) in Theorem~B are weaker, as it suffices to verify the right-hand side of (\ref{eq:RIP}) for $1$-sparse vectors rather than for every $s$-sparse vector.
This happens to be a substantial difference: the assumption
that for every $t\in\Sigma_s$, $\norm{\Gamma t}_2\leq (1+\delta)\norm{t}_2$ is a costly one, and happens to be the reason for the gap between the RIP and the exact reconstruction property. Indeed, while the lower bound in the RIP holds for rather general ensembles (see \cite{Shahar-COLT} and the next section for more details), and is guaranteed solely by the small-ball condition, the upper bound is almost equivalent to having the coordinates of $X$ exhibit a subgaussian behaviour of moments, at least up to some level. Even the fact that one has to verify the upper bound for $1$-sparse vectors comes at a cost, namely, the moment assumption (1) in Theorem~A.

The second goal of this note is to illustrate that while Exact Reconstruction is `cheaper' than the RIP, it still comes at a cost -- namely, that the moment condition (1) in Theorem~A is truly needed.

\begin{Definition}
A random matrix $\Gamma$ is generated by the random variable $x$ if $\Gamma=N^{-1/2}\sum_{i=1}^N \inr{X_i,\cdot}f_i$ and $X_1,...,X_N$ are independent copies of the random vector $X=(x_1,...,x_n)^\top$ whose coordinates are independent copies of $x$.
\end{Definition}
\vskip0.5cm

\noindent {\bf Theorem C.} \textit{ There exist absolute constants
  $c_0,c_1,c_2$ and $c_3$ for which the following holds. Given $n\geq
  c_0$ and $N \log N\leq c_1 n$, there exists a mean-zero, variance
  one random variable $x$ with the following properties:
  \begin{description}
  \item{$\bullet$} $\norm{x}_{L_p} \leq c_2\sqrt{p}$ for $2< p\leq  c_3 (\log n)/(\log N)$.
  \item{$\bullet$} If $(x_j)_{j=1}^n$ are independent copies of $x$
    then $X=(x_1,...,x_n)^\top$ satisfies the small-ball condition
    with constants $u$ and $\beta$ that depend only on $c_2$.
  \item{$\bullet$}  Denote by $\Gamma$ the $N \times n$ matrix generated by $x$. For every $k\in\{1,\ldots,n\}$, with probability larger
    than $1/2$, $\argmin\big(\norm{t}_1:\Gamma t=\Gamma e_k\big)\neq
    \{e_k\}$; therefore, $e_k$ is not exactly reconstructed by Basis
    Pursuit and so $\Gamma$  does not satisfy the exact reconstruction property of order $1$.
  \end{description}
      }

\vskip0.3cm

To put Theorem~C in some perspective, note that if $\Gamma$ is generated by $x$ for which $\norm{x}_{L_2}=1$ and $\norm{x}_{L_p} \leq c_4 \sqrt{p}$ for $2 < p \leq c_5\log  n $, then $X=(x_i)_{i=1}^n$ satisfies the small-ball condition with constants that depend only on $c_4$, and by Theorem~A, if $N \geq c_6 \log n$,  $\Gamma$ satisfies ER($1$) with high probability. On the other hand, the random ensemble from Theorem~C is generated by $x$ that has almost identical properties -- with one exception: its $L_p$ norm is well behaved only for $p \leq c_7 (\log  n)/\log \log n $. This small gap in the number of moments has a significant impact: with probability at least $1/2$, $\Gamma$ does not satisfy ER($1$) when $N$ is of the order of $\log n$.

Therefore, the moment condition in Theorem~A is indeed required (up to a $\log \log n$ factor).
\vskip0.3cm
The idea behind the proof of Theorem~C is to construct a random matrix $\Gamma$ for which, given any basis vector $e_k$, with probability at least $1/2$, $\|\Gamma e_k\|_2 \leq 1$, while the set $\{\Gamma e_j, \ j \not = k\}$ has many `very spiky' vectors: the convex hull ${\rm conv}(\pm \Gamma e_j, \ j \not = k)$ contains a perturbation of $2\sqrt{N}B_1^N$, i.e., a large multiple of the unit ball in $\ell_1^N$. Since such a set must contain the Euclidean unit ball, and in particular, $\Gamma e_k$ as well, it follows that $e_k$ cannot be the unique solution of the $\ell_1$ minimization problem $\min(\|t\|_1 : \ \Gamma t =\Gamma e_k)$. 

The fact that the coordinates of $X$ do not have enough well behaved moments is the key feature that allows one to generate many `spiky' columns in a typical $\Gamma$.

\vskip0.3cm
An alternative formulation of Theorem~C is the following:
\vskip0.3cm
\noindent {\bf Theorem C$^\prime$.} \textit{There are absolute constants $c_0,c_1,c_2$ and $\kappa$ for which the following holds. If $n \geq c_0$ and $2<p<c_1\log n$, there exists a mean-zero and variance $1$ random variable $x$, for which $\|x\|_{L_q} \leq \kappa\sqrt{q}$ for $2 < q \leq p$, and if $N \leq c_2\sqrt{p} (n/\log n)^{1/p}$ and $\Gamma$ is the $N \times n$ matrix generated by $x$, then with probability at least $1/2$, $\Gamma$ does not satisfy the exact reconstruction property of order 1.}
\vskip0.5cm

Theorem~C and Theorem~C$^\prime$ imply that Basis Pursuit may perform poorly when the
coordinates of $X$ do not have enough moments, and requires a
polynomial number of measurements in $n$ to ensure Exact
Reconstruction. This happens to be the price of convex relaxation: a rather striking observation is that $\ell_0$-minimization achieves recovery with the optimal number of measurements under an even weaker small-ball condition than in Theorem~A, and without any additional moment assumptions.

Recall that
$\ell_0$-minimization is defined by  $\hat{x}={\rm argmin} \big(\norm{t}_0:\Gamma
t=\Gamma x_0\big)$, where $\norm{t}_0$ is the cardinality of the
support of $t$.
\begin{Definition} \label{def:weak-small-ball}
$X$ satisfies a weak small-ball condition in $\Sigma_s$ with constant $\beta$
if for every $t\in\Sigma_s$,
\begin{equation}
  \label{eq:weak-small-ball}
 P\big(|\inr{X,t}|>0\big)\geq \beta.
\end{equation}
\end{Definition}

\vskip0.3cm

\noindent {\bf Theorem D.} \textit{For every $0<\beta<1$ there exist constants $c_0$ and $c_1$ that depend only on $\beta$ and for which the following holds. Let $X$ be a random vector that satisfies the weak small-ball condition in $\Sigma_s$ with a constant $\beta$. Let $X_1,\ldots,X_N$ be $N$ independent copies of $X$ and set $\Gamma=N^{-1/2} \sum_{i=1}^N \inr{X_i,\cdot}f_i$. If $N \geq c_0 s \log(en/s)$ then with probability at least $1-2\exp(-c_1N)$, for every $x_0 \in \Sigma_{\lfloor s/2\rfloor}$, $\ell_0$-minimization has a unique solution, which is $x_0$ itself.
}

\vskip0.5cm

The price of convex relaxation can now be clearly seen through the number of measurements needed for exact reconstruction: consider the random vector $X$ constructed in Theorem~C$^\prime$ for, say, $p=4$. Since $X$ satisfies the conditions of Theorem~D, $\ell_0$ minimization may be used to recover any $s$-sparse vector with only $N =cs\log(en/s)$ random measurements. In contrast, Basis Pursuit requires at least $\sim (n/\log n)^{1/4}$ measurements to reconstruct $1$-sparse vectors.

It should be noted that under much stronger assumptions on $X$, the
exact recovery of $s$-sparse vectors using $\ell_0$-minimization may
occur when $N$ is as small as $2s$. Indeed, it suffices to ensure that
all the $N\times 2s$ sub-matrices of $\Gamma$ are non-singular, and
this is the case when $N=2s$ if the entries of $\Gamma$ are
independent random variables that are absolutely continuous (see
Chapter~2 in \cite{MR3100033} for more details).

We end this introduction with a word about notation and the
organization of the article. The proofs of Theorem~A, Theorem~B and
Theorem~D are presented in the next section, while the proofs of
Theorem~C and Theorem~C$^\prime$ may be found in
Section~\ref{sec:proof_theo_C}. The final section is devoted to
results in a natural `noisy' extension of Compressed Sensing. In
particular, we prove that both the {\it Compatibility Condition} and
the {\it Restricted Eigenvalue Condition} hold under weak moment
assumptions; we also study related properties of random polytopes.

As for notation, throughout, absolute constants or constants that depend on other parameters are denoted by $c$, $C$, $c_1$,
$c_2$, etc., (and, of course, we will specify when a constant is
absolute and when it depends on other parameters). The values of these
constants may change from line to line. The notation $x\sim y$ (resp. $x\lesssim y$) means that there exist absolute constants $0<c<C$ for which $cy\leq x\leq Cy$ (resp. $x\leq Cy$). If $b>0$ is a parameter then $x\lesssim_b y$ means that $x\leq C(b) y $ for some constant $C(b)$ that depends only on $b$.

Let $\ell_p^m$ be $\R^m$ endowed with the norm $\|x\|_{\ell_p^m}=\big(\sum_{j}|x_j|^p\big)^{1/p}$; the corresponding unit ball is denoted by $B_p^m$ and the unit Euclidean sphere in $\R^m$ is $S^{m-1}$. If $A\subset\R^n$ then $\IND_A$ denotes the indicator function of $A$. Finally, we will assume that $(\cX,\mu)$ is a probability space, and that $X$ is distributed according to $\mu$.

\section{Proof of Theorem~A, B and D}
\label{sec:proof-thm-A}
The proof of Theorem~A has several components, and although the first of which is rather standard, we present it for the sake of completeness.
\begin{Lemma}\label{prop:exact-reconstruct}
Let $\Gamma:\R^n \to \R^N$ be a matrix and set ${\rm ker}(\Gamma)$ to be its kernel. If $0<r<1$ and $B_1^n \cap r S^{n-1}$ does not intersect ${\rm ker}(\Gamma)$, then $\Gamma$ satisfies the exact reconstruction property in $\Sigma_{\lfloor (2r)^{-2}\rfloor}$.
\end{Lemma}
\proof Observe that if $x \in B_1^n$ and $\|x\|_2  \geq r$ then $y=rx/\|x\|_2 \in B_1^n \cap rS^{n-1}$. Therefore, if $y \not \in {\rm ker}(\Gamma)$, the same holds for $x$; thus
$$
\sup_{x \in B_1^n \cap {\rm ker}(\Gamma)} \|x\|_2 <r.
$$
Let $s=\lfloor (2r)^{-2}\rfloor$, fix $x_0 \in \Sigma_s$ and put $I$ to be the set indices of coordinates on which $x_0$ is supported. Given a nonzero $h \in {\rm ker}(\Gamma)$, let $h=h_I+h_{I^c}$ -- the decomposition of $h$ to coordinates in $I$ and in $I^c$. Since $h/\|h\|_1 \in B_1^n \cap {\rm ker}(\Gamma)$, it follows that $\|h\|_2 < r \|h\|_1$, and by the choice of $s$, $2\sqrt{s}\|h\|_2 < \|h\|_1$.
Therefore,
\begin{align*}
  \|x_0+h\|_1&= \|x_0+h_I\|_1+\|h_{I^c}\|_1\geq \|x_0\|_1-\|h_I\|_1+\|h_{I^c}\|_1
  \\
  & = \|x_0\|_1-2\|h_I\|_1+\|h\|_1 \geq \|x_0\|_1-2\sqrt{|I|}\|h_I\|_2+\|h\|_1 > \|x_0\|_1.
\end{align*}
Hence, $\|x_0+h\|_1>\|x_0\|_1$ and $x_0$ is the unique minimizer of the basis pursuit algorithm.
\endproof

The main ingredient in the proof of Theorem~A is Lemma \ref{lemma:small-ball} below, which is based on the small-ball method introduced in \cite{Shahar-COLT,Shahar-Gelfand}. To formulate the lemma, one requires the notion of a VC class of sets.

\begin{Definition}
Let $\cG$ be a class of $\{0,1\}$-valued functions defined on a set $\cX$. The set $\cG$ is a VC-class if there exists an integer $V$ for which, given any $x_1,...,x_{V+1}\in\cX$,
\begin{equation} \label{eq:vc}
\left|\left\{(g(x_1),...,g(x_{V+1})):g\in\cG\right\}\right|<2^{V+1}.
\end{equation}
The VC-dimension of $\cG$, denoted by $VC(\cG)$, is the smallest integer $V$ for which \eqref{eq:vc} holds.
\end{Definition}

The VC dimension is a combinatorial complexity measure that may be used to control the $L_2(\mu)$-covering numbers of the class; indeed, set $N(\cG,\eps,L_2(\mu))$ to be the smallest number of open balls of radius $\eps$ relative to the $L_2(\mu)$ norm that are needed to cover $\cG$. A well known result due to Dudley \cite{MR512411} is that if $VC(\cG)=V$ and $\mu$ is a probability measure on $\cX$ then for every $0<\eps<1$,
\begin{equation}
  \label{eq:VC-entropy}
N\big(\cG,\eps,L_2(\mu)\big)\leq \Big(\frac{c_1}{\eps}\Big)^{c_2V},
\end{equation}
where $c_1$ and $c_2$ are absolute constants.

\begin{Lemma} \label{lemma:small-ball}
There exist absolute constants $c_1$ and $c_2$ for which the following holds.
Let $\cF$ be a class of functions and assume that there are $\beta>0$ and $u \geq 0$ for which
$$
\inf_{f \in \cF} P\big(|f(X)|>u\big) \geq \beta.
$$
Let $\cG_u = \left\{ \IND_{\{|f|>u\}} : f \in \cF\right\}$. If $VC(\cG_u) \leq d$ and $N \geq c_1 d/\beta^2$ then with probability at least $1-\exp(-c_2\beta^2 N)$,
$$
\inf_{f \in \cF} \big|\big\{ i\in\{1,\ldots,N\}: |f(X_i)| > u\big\}\big| \geq \frac{\beta N}{2}.
$$
\end{Lemma}
\begin{Remark}
Note that $u=0$ is a `legal choice' in Lemma \ref{lemma:small-ball}, a fact that will be used in the proof of Theorem D.
\end{Remark}

\proof Let $G(X_1,...,X_N)=\sup_{g \in \cG_u} |N^{-1}\sum_{i=1}^N g(X_i)-\E g(X)|$. By the bounded differences inequality (see, for example, Theorem~6.2 in \cite{BouLugMass13}), with probability at least $1-\exp(-t)$,
$$
G(X_1,...,X_N) \leq \E G(X_1,...,X_N) + c_1 \sqrt{\frac{t}{N}}.
$$
Standard empirical processes arguments (symmetrization, the fact that Bernoulli processes are
subgaussian and the entropy estimate \eqref{eq:VC-entropy} -- see, for example, Chapters~2.2, 2.3 and 2.6 in \cite{vanderVaartWellner}), show that since $VC(\cG) \leq d$,
\begin{equation}\label{eq:vc_bound_dudley}
\E G(X_1,...,X_N) \leq c_2\sqrt{\frac{d}{N}} \leq \frac{\beta}{4},
\end{equation}
provided that $N \gtrsim d/\beta^2$. Therefore, taking $t=N\beta^2/16c_1^2$, it follows that with probability at least $1-\exp(-c_3\beta^2 N)$, for every $f \in \cF$,
$$
\frac{1}{N} \sum_{i=1}^N \IND_{\{|f|>u\}}(X_i) \geq P\big(|f(X)|>u\big)-\frac{\beta}{2} \geq \frac{\beta}{2}.
$$
Therefore, on that event, $|\{i : |f(X_i)|>u\}| \geq \beta N/2$ for every $f \in \cF$.
\endproof

\begin{Corollary} \label{cor:matrix}
There exist absolute constants $c_1$ and $c_2$ for which the following holds.
Let $X \in \R^n$ be a random vector.
\begin{description}
\item{1.} If there are $0<\beta\leq1$ and $u\geq0$ for which
$P\big(|\inr{t,X}|>u\big) \geq \beta$ for every $t \in S^{n-1}$ and if $N \geq c_1n/\beta^2$, then with probability at least $1-\exp(-c_2N\beta^2)$,
$$
\inf_{t \in S^{n-1}} \frac{1}{N} \sum_{i=1}^N \inr{X_i,t}^2 > \frac{u^2 \beta}{2}.
$$
\item{2.} If there are $0<\beta\leq1$ and $u\geq0$ for which $P\big(|\inr{t,X}|>u\big) \geq \beta$ for every $t \in \Sigma_s\cap S^{n-1} $ and if $N \geq c_1s\log(en/s)/\beta^2$, then with probability at least $1-\exp(-c_2N\beta^2)$,
$$
\inf_{t \in \Sigma_s\cap S^{n-1}} \frac{1}{N} \sum_{i=1}^N \inr{X_i,t}^2 > \frac{u^2 \beta}{2}.
$$
\end{description}
\end{Corollary}

\begin{Remark}
Note that the first part of Corollary \ref{cor:matrix} gives an estimate on the smallest singular value of the random matrix $\Gamma=N^{-1/2}\sum_{i=1}^N \inr{X_i,\cdot}f_i$. The proof follows the same path as in \cite{Shahar-Vladimir}, but unlike the latter, no assumption on the covariance structure of $X$, used both in \cite{Shahar-Vladimir} and in \cite{MR3127875}, is required. In fact, Corollary \ref{cor:matrix} may be applied even if the covariance matrix does not exist.  Thus, under a small-ball condition, the smallest singular value of $\Gamma$ is larger than $c(\beta,u)$ with high (exponential) probability.
\end{Remark}

\noindent{\bf Proof of Corollary \ref{cor:matrix}.} To prove the first part of the claim, let $\cF=\{\inr{t,\cdot} : t \in S^{n-1}\}$. Recall that the VC dimension of a class of half-spaces in $\R^n$ is at most $n$, and thus, one may verify that for every $u \geq 0$, the VC dimension of $$
\cG_u=\{\IND_{\{|f| > u\}} : f \in \cF \}
$$
is at most $c_1n$ for a suitable absolute constant $c_1$ (see, e.g., Chapter~2.6 in \cite{vanderVaartWellner}). The claim now follows immediately from Lemma~\ref{lemma:small-ball} because $$
\frac{1}{N} \sum_{i=1}^N \inr{t,X_i}^2 > \frac{u^2}{N} |\{i: |\inr{X_i,t}| > u\}|
$$
for every $t\in S^{n-1}$.

Turning to the second part, note that $\Sigma_s \cap S^{n-1}$ is a union of $\binom{n}{s}$ spheres of dimension $s$. Applying the first part to each one of those spheres, combined with the union bound, it follows that for $N \geq c_2 \beta^{-2}s\log(en/s)$, with probability at least $1-\exp(-c_3N\beta^2)$,
$$
\inf_{t \in \Sigma_s \cap S^{n-1}} \frac{1}{N}\sum_{i=1}^N \inr{X_i,t}^2 > \frac{u^2 \beta}{2}.
$$\endproof

Corollary \ref{cor:matrix} shows that the small-ball condition for linear functionals implies that $\Gamma$ `acts well' on $s$-sparse vectors. However, according to Lemma~\ref{prop:exact-reconstruct}, exact recovery is possible if $\Gamma$ is well behaved on the set
$$
\sqrt{\kappa_0 s} B_1^n \cap S^{n-1} = \{x \in \R^n : \|x\|_1 \leq \sqrt{\kappa_0s}, \ \ \|x\|_2 =1\}
$$
for a well-chosen constant $\kappa_0$. In the standard (RIP-based) argument, one proves exact reconstruction by first showing that the RIP holds in $\Sigma_s$, and then the fact that each vector in $\sqrt{\kappa_0 s} B_1^n \cap S^{n-1}$ is well approximated by vectors from $\Sigma_s$ (see, for instance, \cite{MR3113826}) allows one to extend the RIP from $\Sigma_s$ to $\sqrt{\kappa_0 s} B_1^n \cap S^{n-1}$. Unfortunately, this extension requires {\it both upper and lower} estimates in the RIP.

Since the upper part of the RIP in $\Sigma_s$ forces severe restrictions on the random vector $X$, one has to resort to a totally different argument if one wishes to extend the {\it lower bound} from $\Sigma_s$ (which only requires the small-ball condition) to $\sqrt{\kappa_0 s} B_1^n \cap S^{n-1}$.

The method presented below is based on Maurey's empirical method and has been recently used in \cite{Oliveira13}.
\begin{Lemma} \label{lemma:Maurey}
Let $\Gamma:\R^n \to \R^N$, put $1<s\leq n$ and assume that for every $x \in \Sigma_s$, $\|\Gamma x\|_2 \geq \lambda \|x\|_2$. If $y \in \R^n$ is a nonzero vector and $\mu_j = |y_j|/\|y\|_1$, then
$$
\|\Gamma  y\|_2^2 \geq \lambda^2 \|y\|_2^2 - \frac{\|y\|_1^2}{s-1} \left( \sum_{j=1}^n \norm{\Gamma e_j}_2^2 \mu_j- \lambda^2\right).
$$
\end{Lemma}

\proof
Fix $y \in \R^n$, let $Y$ be a random vector in $\R^n$ defined by
$$
P(Y=\|y\|_1\sgn(y_j)e_j)=|y_j|/\|y\|_1,
$$ for every $j=1,\ldots,n$ and observe that $\E Y =y$.

Let $Y_1,...,Y_s$ be independent copies of $Y$ and set $Z=s^{-1}\sum_{k=1}^s Y_k$. Note that $Z \in \Sigma_s$ for every realization of $Y_1,...,Y_s$; thus
$\|\Gamma Z\|_2^2 \geq \lambda^2 \|Z\|_2^2$ and
\begin{equation} \label{eq:lower-Maurey}
\E \|\Gamma Z\|_2^2 \geq \lambda^2 \E\|Z\|_2^2.
\end{equation}
It is straightforward to verify that $\E \inr{Y,Y} = \|y\|_1^2$; that  if $i \not = j$ then $\E \inr{\Gamma Y_i,\Gamma Y_j} = \inr{\Gamma y,\Gamma y}$; and that for every $1 \leq k \leq s$,
$$
\E \inr{\Gamma Y_k,\Gamma Y_k} = \|y\|_1 \sum_{j=1}^n |y_j| \norm{\Gamma e_j}_2^2.$$

Therefore, setting $\mu_j = |y_j|/\|y\|_1$ and $W = \sum_{j=1}^n \norm{\Gamma e_j}_2^2\mu_j$,
\begin{align*}
\E \|\Gamma Z\|_2^2 = & \frac{1}{s^2} \sum_{i,j=1}^s \E \inr{\Gamma Y_i,\Gamma Y_j} = \left(1-\frac{1}{s}\right) \|\Gamma y\|_2^2 + \frac{\|y\|_1}{s} \sum_{j=1}^n |y_j| \norm{\Gamma e_j}_2^2
\\
= & \left(1-\frac{1}{s}\right) \|\Gamma y\|_2^2 +  W\frac{\|y\|_1^2}{s},
\end{align*}
and using the same argument one may show that
$$
\E\|Z\|_2^2 = \left(1-\frac{1}{s}\right)\|y\|_2^2 +\frac{\|y\|_1^2}{s}.
$$
Combining these two estimates with \eqref{eq:lower-Maurey},
$$
\left(1-\frac{1}{s}\right) \|\Gamma y\|_2^2 \geq \lambda^2 \left(\left(1-\frac{1}{s}\right)\|y\|_2^2 +\frac{\|y\|_1^2}{s}\right) -  W\frac{\|y\|_1^2}{s},
$$
proving the claim.
\endproof

\textbf{Proof of Theorem~B:} Assume that for every $x \in \Sigma_s$, $\|\Gamma x\|_2 \geq c_0 \|x\|_2$ and that for every $1 \leq i \leq n$, $\|\Gamma e_i\|_2 \leq c_1$. It follows from Lemma~\ref{lemma:Maurey} that if $s-1>c_1^2/(c_0^2 r^2)$, then for every $y\in B_1^n\cap r S^{n-1}$,
\begin{equation*}
  \norm{\Gamma y}_2^2\geq c_0^2 \norm{y}_2^2-\frac{\norm{y}_1}{s-1}\sum_{i=1}^n \norm{\Gamma e_i}_2^2|y_i|\geq c_0^2 r^2 -\frac{c_1^2}{s-1}>0.
\end{equation*}
The claim now follows from Lemma~\ref{prop:exact-reconstruct}.
\endproof

Consider the matrix $\Gamma=N^{-1/2}\sum_{i=1}^N\inr{X_i,\cdot}f_i$. Observe that for every $t \in \R^n$, $\|\Gamma t\|_2^2= N^{-1}\sum_{i=1}^N \inr{X_i,t}^2$, and that if $X_j=(x_{i,j})_{i=1}^n$ then
$$
\norm{\Gamma e_j}_2^2= \frac{1}{N}\sum_{i=1}^N x_{i,j}^2,
$$
which is an average of $N$ iid random variables (though $\norm{\Gamma e_1}_2,$ $\ldots,\norm{\Gamma e_n}_2$ need not be independent).

Thanks to Theorem~B and Corollary \ref{cor:matrix}, the final component needed for the proof of Theorem~A is information on the sum of iid random variables, which will be used to bound $\max_{1\leq j\leq n}\norm{\Gamma e_j}_2^2$ from above.

\begin{Lemma}\label{prop:moment-sum-variid}
 There exists an absolute constant $c_0$ for which the following holds.
  Let $z$ be a mean-zero random variable and put $z_1,\ldots,z_N$ to
  be $N$ independent copies of $z$. Let $p_0\geq 2$ and assume that
  there exists $\kappa_1>0$ and $\alpha \geq 1/2$ for which
  $\|z\|_{L_p}\leq \kappa_1 p^\alpha$ for every $2\leq p\leq p_0$. If
  $N\geq p_0^{\max\{2\alpha-1,1\}}$  then for every $2\leq p\leq p_0$,
    \begin{equation*}
      \Big\|\frac{1}{\sqrt{N}}\sum_{i=1}^N z_i\Big\|_{L_p}\leq c_1(\alpha) \kappa_1 \sqrt{p},
    \end{equation*}
  where $c_1(\alpha)=c_0\exp((2\alpha-1))$.
\end{Lemma}

Lemma~\ref{prop:moment-sum-variid} shows that even under a  weak moment assumption, namely that $\|z\|_{L_p} \lesssim p^{\alpha}$ for $p \leq p_0$ and $\alpha \geq 1/2$ that can be large, a normalized sum of $N$ independent copies of $z$ exhibits a `subgaussian' moment growth up to the same $p_0$, as long as $N$ is sufficiently large.

The proof of Proposition~\ref{prop:moment-sum-variid}  is based on the following fact due to Lata{\l}a.
\begin{Theorem}[\cite{MR1457628}, Theorem~2 and Remark~2]\label{thm:latala}
If $z$ is a mean-zero random variable and $z_1,...,z_N$ are independent copies of $z$, then for any $p\geq2$,
$$
\big\|\sum_{i=1}^N z_i\big\|_{L_p} \sim \sup \left\{\frac{p}{s}\left(\frac{N}{p}\right)^{1/s}\|z\|_{L_s} : \max\{2,p/N\} \leq s \leq p \right\}.
$$
\end{Theorem}

\noindent{\bf Proof of Lemma~\ref{prop:moment-sum-variid}.} Let $2\leq
p\leq p_0$ and $N\geq p$. Since $\|z\|_{L_s} \leq \kappa_1 s^{\alpha}$
for any $2\leq s\leq p$, it follows from Theorem \ref{thm:latala} that
\begin{equation*}
\big\|\sum_{i=1}^N z_i\big\|_{L_p} \leq c_0 \kappa_1 \sup\left\{ p(N/p)^{1/s} s^{-1+\alpha}:\max\{2,p/N\} \leq s \leq p\right\}.
\end{equation*}
It is straightforward to verify that the function $h(s)=(N/p)^{1/s}
s^{-1+\alpha}$ is non-increasing when $\alpha \leq 1$ and attains its
maximum in $s=\max\{2,p/N\}=2$ or in $s=p$ when $\alpha>1$. Therefore, when $N \geq p$,
$$
\big\|\sum_{i=1}^N z_i\big\|_{L_p} \leq c_1 \kappa_1 \max\left\{\sqrt{N p}, N^{1/p} p^{\alpha} \right\}.
$$
Finally, if $N \geq p^{2\alpha-1}$ then $e^{2\alpha-1} \sqrt{Np} \geq N^{1/p} p^{\alpha}$, which completes the proof.
\endproof

\noindent{\bf Proof of Theorem A.}  Consider $N \geq c_1s\log(en/s)/\beta^2$. By Corollary \ref{cor:matrix}, with probability at least $1-\exp(-c_2N\beta^2)$,
\begin{equation} \label{eq:in-proof-A-1}
\inf_{t \in \Sigma_s\cap S^{n-1}} \frac{1}{N} \sum_{i=1}^N \inr{X_i,t}^2 > \frac{u^2 \beta}{2}.
\end{equation}
Set $(X_i)_{i=1}^N$ for which \eqref{eq:in-proof-A-1} holds and let $\Gamma = N^{-1/2}\sum_{i=1}^N \inr{X_i,\cdot}f_i$. By Lemma \ref{lemma:Maurey} for $\lambda^2=u^2 \beta/2$, it follows that when $r \geq 1$,
\begin{equation} \label{eq:in-proof-Maurey}
\inf_{t \in \sqrt{r}B_1^n \cap S^{n-1}} \|\Gamma t\|_2^2 \geq \lambda^2-\frac{2r}{s} \max_{1 \leq j \leq n} \norm{\Gamma e_j}_2^2.
\end{equation}
Next, one has to obtain a high probability upper estimate on $\max_{1 \leq j \leq n} \|\Gamma e_j\|_2^2$. To that end, fix $w \geq 1$ and consider $z=x_j^2-1$ - where $x_j$ is the
$j$-th coordinate of $X$. Observe that $z$ is a centered random variable and that $\norm{z}_{L_q}\lesssim 4^\alpha\kappa_1^2 q^{2\alpha}$ for every $1\leq q\leq \kappa_2\log(wn)$. Thus, by Lemma \ref{prop:moment-sum-variid} for $p = \kappa_2\log(wn)$ and
$c_3(\alpha) \sim 4^\alpha\exp((4\alpha-1))$,
\begin{equation*}
\big\|\frac{1}{N}\sum_{i=1}^N z_i\big\|_{L_p}\leq c_3(\alpha) \kappa_1^2 \sqrt{\frac{p}{N}},
\end{equation*}
provided that $ N\geq p^{\max\{4\alpha-1,1\}}=(\kappa_2\log(wn))^{\max\{4\alpha-1,1\}}$.
Hence, if $N \geq (c_3(\alpha)\kappa_1^2)^2 (\kappa_2\log(wn))^{\max\{4\alpha-1,1\}}$, and setting  $V_j=\norm{\Gamma e_j}_2^2$, one has
$$
\|V_j\|_{L_p}=\|\frac{1}{N}\sum_{i=1}^N x_{i,j}^2 \|_{L_p} \leq 1+c_3(\alpha) \kappa_1^2 \sqrt{\frac{\kappa_2\log(wn)}{N}} \leq 2;
$$
thus,
\begin{align*}
P(\max_{1 \leq j \leq n} V_j \geq 2e) \leq & \sum_{j=1}^n P(V_j \geq 2e) \leq  \sum_{j=1}^n \left(\frac{\|V_j\|_{L_p}}{2e}\right)^p
\\
\leq & n\left(\frac{1}{e}\right)^p=\frac{1}{w^{\kappa_2}n^{\kappa_2-1}}.
\end{align*}
Combining the two estimates, if
$$
N \gtrsim \max\left\{s \log (en/s), (c_3(\alpha)\kappa_1^2)^2 (\kappa_2\log(wn))^{\max\{4\alpha-1,1\}} \right\}
$$
and $r \leq s \lambda^2/8e=s u^2 \beta/16e$, then with probability at least $1-\exp(-c_2N\beta^2)-1/(w^{\kappa_2} n^{\kappa_2-1})$,
\begin{equation} \label{eq:in-proof-cone}
\inf_{t \in \sqrt{r}B_1^n \cap S^{n-1}} \|\Gamma t\|_2^2 \geq \lambda^2-\frac{4e r}{s} \geq \lambda^2/2.
\end{equation}
Therefore, by Lemma \ref{prop:exact-reconstruct}, $\Gamma$ satisfies the exact reconstruction property for vectors that are $c_4u^2\beta s$-sparse, as claimed.
\endproof


\textbf{Proof of Theorem~D:} Since the argument is almost identical to the one used in the proof of the second part of Corollary \ref{cor:matrix}, we will only sketch the details. Observe that if $\Gamma=N^{-1/2}\sum_{i=1}^N \inr{X_i,\cdot}f_i$ and ${\rm ker}(\Gamma) \cap \Sigma_s=\{0\}$, then for any $x_0\in\Sigma_{\lfloor s/2\rfloor}$, the only $t\in\R^n$ for which $\Gamma t=\Gamma x_0 $ and $\norm{t}_0\leq \norm{x_0}_0$ is $x_0$ itself. Thus, it suffices to show that for every $x \in \Sigma_s \cap S^{n-1}$, $|\inr{X_i,x}|>0$ for some $1 \leq i \leq n$. Since $\Sigma_s \cap S^{n-1}$ is a union of $\binom{n}{s}$ spheres, the claim follows from Lemma~\ref{lemma:small-ball} applied to each one of those spheres and for $u=0$, combined with a union bound argument.
\endproof

\section{Proof of Theorem~C and Theorem~C$^\prime$}
\label{sec:proof_theo_C}
Consider an $N \times n$ matrix $\Gamma$ and $J\subset\{1,\ldots,n\}$. Set $\Gamma_J$
to be the  $(N\times|J|)$ restriction of $\Gamma$ to ${\rm span}\{e_j : j \in J\}$. Recall that $B_1^n$ is the unit ball in $\ell_1^n=(\R^n, \| \cdot \|_1)$, and put $B_1^{J^c}$ to be the set of vectors in $B_1^n$ that are supported in $J^c$ -- the complement of $J$ in $\{1,...,n\}$.

\begin{Lemma}\label{lem:non_ER}
Fix integers $s,N \leq n$ and $J \subset \{1,...,n\}$ of cardinality at most $s$. If $v\in\R^n$ is supported in $J$, $\|v\|_1=1$ and $\Gamma v \in\Gamma B_1^{J^c}$, then  $\Gamma$ does not satisfy the exact reconstruction property of order $s$.
\end{Lemma}
\proof
Let $w\in B_1^{J^c}$ for which $\Gamma v=\Gamma w$ and observe that $v \not = w$ (otherwise, $v\in B_1^J\cap B_1^{J^c}$, implying that $v=0$, which is impossible because $\norm{v}_1=1$).

Since $\|w\|_1 \leq 1 =\|v\|_1$, $w$ is at least as good a candidate as $v$ for the $\ell_1$-minimization problem $\min\big(\norm{t}_1:\Gamma t=\Gamma v\big)$; hence, $v$ is not the unique solution of that problem.
\endproof

Set $x_{\cdot1},\cdots,x_{\cdot n}$ to be the columns of $\Gamma$.
It immediately follows from Lemma~\ref{lem:non_ER} that if one wishes
to prove that $\Gamma$ does not satisfy ER($1$), it suffices to show
that, for instance the first basis vector $e_1$ cannot be exactly
reconstruct. This follows from
\begin{equation*}
 \Gamma e_1 = x_{\cdot 1} \in {\rm absconv}\big( x_{\cdot k}:k\neq
 1\big) = {\rm absconv}\big(\Gamma e_k: k \not = 1 \big)=\Gamma B_1^{\{1\}^c},
\end{equation*}
where $ {\rm absconv}(S)$ is the convex hull of $S \cup -S$. Therefore, if
\begin{equation}
  \label{eq:er_1_4}
  \norm{x_{\cdot 1}}_2\leq c_0 \ \ \mbox{ and } \ \  c_0B_2^N \subset {\rm absconv}\big( x_{\cdot k}:k\neq 1\big),
\end{equation} for some absolute constant $c_0$,  then $\Gamma$ does not satisfy ER($1$).

The proofs of Theorem C and of Theorem C$^\prime$ follow from the
construction of a random matrix ensemble for which (\ref{eq:er_1_4})
holds with probability larger than $1/2$. We now turn on to such a construction.

Let $\eta$ be a selector (a $\{0,1\}$-valued random variable) with mean $\delta$ to be named later, and let $\eps$ be a symmetric $\{-1,1\}$-valued random variable that is independent of $\eta$. Fix $R>0$ and set
$$
z= \eps(1+R  \eta).
$$

Observe that if $p\geq2$ and $R\geq 1$ then
\begin{equation*}
  \frac{\norm{z}_{L_p}}{\norm{z}_{L_2}}= \frac{\big(1+ \big((1+R)^p-1\big)\delta\big)^{1/p}}{\big(1+ \big((1+R)^2-1\big)\delta\big)^{1/2}}\sim \frac{(1+R^p\delta)^{1/p}}{(1+R^2\delta)^{1/2}}\sim R \delta^{1/p},
\end{equation*}
and the last equivalence holds when $R^2\delta\lesssim 1$ and  $R^p\delta\gtrsim 1$. Fix $2<p \leq 2\log(1/\delta)$ which will be specified later and set $R=\sqrt{p}(1/\delta)^{1/p}$. Since the function $q \to \sqrt{q}/\delta^{1/q}$ is decreasing for $2 \leq q \leq 2\log(1/\delta)$ one has that for $2 \leq q \leq p$ and for $\delta$ that is small enough,
$$
\|z\|_{L_q} \leq c_0\sqrt{q}\|z\|_{L_2}.
$$

Note that $x=z/\norm{z}_{L_2}$ is a mean-zero, variance one random
variable that exhibits a `subgaussian' moment behaviour only up to
$p$. Indeed, if $2 \leq q \leq p$, $\|z\|_{L_q} \lesssim
\sqrt{q}\|z\|_{L_2}$, and if $q> p$, $\norm{z}_{L_q}\sim
\sqrt{p}\delta^{1/q-1/p}\norm{z}_{L_2}$, which may be far larger
than $\sqrt{q}\norm{z}_{L_2}$ if $\delta$ is sufficiently small.

Let $X=(x_1,\ldots,x_n)$ be a vector whose coordinates are independent, distributed as $x$ and let $\Gamma$ be the measurement matrix generated by $x$.
Note that up to the normalization factor of $\norm{z}_{L_2}$, which is of the order of a constant when $R^2 \delta\lesssim 1$, $\sqrt{N}\Gamma$ is a perturbation of a Rademacher matrix by a sparse matrix with few random spikes that are either $R$ or $-R$.

As noted earlier, if for every $t \in \R^n$,
\begin{equation} \label{eq:L-4-L-2}
\|\inr{X,t}\|_{L_4} \leq C\|\inr{X,t}\|_{L_2},
\end{equation}
then the small-ball condition holds with constants that depend only on $C$. To show that $X$ satisfies \eqref{eq:L-4-L-2}, denote by $\E_\eta$ (resp. $\E_\eps$) the expectation with respect to the $\eta$-variables (resp. $\eps$-variables), and observe that by a straightforward application of Khintchine's inequality (see, e.g., p.91 in \cite{LT:91}), for every $t\in\R^n$,
\begin{align*}
  &\E\inr{X,t}^4\lesssim \E_\eta \E_\eps\Big(\sum_{j=1}^n \eps_j(1+R\eta_j)t_j\Big)^4\lesssim \E_\eta\Big(\sum_{j=1}^n (1+R\eta_j)^2t_j^2 \Big)^2
  \\
&=\E_\eta \sum_{k,\ell}(1+R\eta_k)^2t_k^2(1+R\eta_\ell)^2t_\ell^2\lesssim \norm{t}_2^4=\Big(\E\inr{X,t}^2\Big)^2,
\end{align*}
provided that $R^4\delta \lesssim 1$.

Set $(f_i)_{i=1}^N$ to be the canonical basis of $\R^N$ and put $\tilde{\Gamma}=(z_{\ell k})=\norm{z}_{L_2}\sqrt{N}\Gamma$, an $N \times n$ matrix whose entries are independent copies of $z$. Let
$$
v_j=\tilde{\Gamma}e_j = \sum_{\ell=1}^N z_{\ell j}f_j,
$$
and consider
$$
V={\rm absconv}\left(\{v_j \ : \ 2 \leq j \leq n\}\right),
$$
the convex hull of $(\pm v_j)_{j=2}^n$.

We will show that with probability at least $1/2$, $\sqrt{N} B_2^N \subset V$ and $\|v_1\|_2 \leq \sqrt{N}$, in three steps:
\begin{Lemma} \label{lemma:perturbation}
With probability at least $3/4$, for every $1 \leq i \leq N$ there is $y_i \in B_\infty^N$ for which $y_i + R f_i \in V$.
\end{Lemma}
In other words, with non-trivial probability, $V$ contains a
perturbation of all the vectors $Rf_i, i=1,\ldots,N$, and thus, $V$ `almost' contains $RB_1^N$.
\vskip0.3cm
\proof Fix a realization of the $N \times n$ Rademacher matrix $(\eps_{\ell j})$ and note that for every $1 \leq i \leq N$ and every $2 \leq j \leq n$
$$
v_j = \sum_{\ell=1}^N \eps_{\ell j} f_\ell + \eps_{ij}Rf_i
$$
if $\eta_{ij}=1$ and for every $\ell \not = i$, $\eta_{\ell j}=0$. Moreover, if this happens, and since $V$ is centrally symmetric (that is, if $v \in V$ then $-v \in V$),
$$
\eps_{ij}\left(\sum_{\ell=1}^N \eps_{\ell j} f_\ell\right)+ Rf_i = y_i+Rf_i \in V,
$$
and $y_i \in B_\infty^N$.

Thus, it remains to estimate the probability that for every $1 \leq i \leq N$ there is some $2 \leq j \leq n$ for which $\eta_{ij}=1$ and for every $\ell \not = i$, $\eta_{\ell j}=0$. Clearly, for every $1 \leq i \leq N$,
\begin{align*}
\  & P_\eta\Big(\mbox{ there exists } j\in\{2,\ldots,n\}: \eta_{ij}=1, \ {\rm and} \ \eta_{\ell j}=0 \ {\rm if} \ \ell \not = i \Big)
\\
&=1-(1-(1-\delta)^{N-1}\delta)^{n-1}\geq 1-\frac{1}{4N}
\end{align*}
provided that
\begin{equation*}
  \frac{\log N}{n}\lesssim \delta \lesssim \frac{\log\big(en/N\big)}{N}.
\end{equation*}
Hence, the claim follows by the union bound and integration with respect to the
$(\eps_{ij})$.
\endproof
Next, it is straightforward to verify that when $V$ contains such a
perturbation of $RB_1^N$ (by vectors in $B_\infty^N$), it must also contain a large Euclidean ball, assuming that $R$ is large enough.

\begin{Lemma} \label{lemma-width}
Let $R>N$, and for every $1 \leq i \leq N$, set $y_i \in B_\infty^N$ and put $v_i=Rf_i+y_{i}$. If $V$ is a convex, centrally symmetric set, and if $v_i \in V$ for every $1 \leq i \leq N$ then $\big(R/\sqrt{N}-\sqrt{N}\big) B_2^N \subset V$.
\end{Lemma}

\proof
A separation argument shows that if $\sup_{v \in V} |\inr{v,w}| \geq \rho$ for every $w \in S^{N-1}$, then $\rho B_2^N \subset V$ (indeed, otherwise there would be some $x \in \rho B_2^N \backslash V$; but it is impossible to separate $x$ and the convex and centrally symmetric $V$ using any norm-one functional).

To complete the proof, observe that for every $w \in S^{N-1}$,
\begin{align*}
&\sup_{v \in V} |\inr{v,w}| \geq  \max_{1 \leq i \leq N} |\inr{Rf_i+y_{i},w}| \\
&\geq \max_{1 \leq i\leq N} |\inr{Rf_i,w}| - \max_{1 \leq i \leq N}|\inr{y_i,w}|
\geq  R/\sqrt{N} - \sqrt{N}.
\end{align*}
\endproof

Applying Lemma~\ref{lemma-width}, it follows that if $R\geq 2N$ then with probability at least $3/4$,
$ \sqrt{N} B_2^N \subset V$. Finally, if $\delta \lesssim 1/N$ then
\begin{equation*}
P\big(\sum_{\ell=1}^N z_{\ell 1} f_\ell \in \sqrt{N}B_2^N\big) \geq P\big(\|\sum_{\ell=1}^N z_{\ell 1} f_\ell \|_2=\sqrt{N}\big)=(1-\delta)^N \geq 3/4.
\end{equation*}
Hence, with probability at least $1/2$,
$$
\sum_{\ell=1}^N z_{\ell 1} f_\ell =\tilde{\Gamma} e_1 \in V= {\rm absconv}\big(\tilde{\Gamma}e_j : j \in \{2,...,n\}\big),
$$
and the same assertion holds for the normalized matrix $\Gamma$, showing that it does not satisfy ER($1$).

Of course, this assertion holds under several conditions on the parameters involved: namely, that $R=\sqrt{p}(1/\delta)^{1/p}\geq 2N$; that $(\log N)/n\lesssim \delta \lesssim\log\big(en/N\big)/N$; that $R^4\delta\lesssim 1$; that $p \leq 2\log(1/\delta)$ and that $\delta\lesssim 1/N$.

For instance, one may select $\delta \sim (\log N)/n$ and $p\sim (\log n)/\log N$, in which case all these conditions are met; hence, with probability at least $1/2$, $\Gamma$ does not satisfy ER($1$), proving Theorem C. A similar calculation leads to the proof of Theorem C$^\prime$.
\endproof
\begin{Remark}
Note that the construction leads to a stronger, non-uniform result,
namely, that for every basis vector $e_k$, with probability at least
$1/2$, $e_k$ is not the unique solution of $\min(\|t\|_1 : \ \Gamma t
= \Gamma e_k)$. In particular, uniformity over all supports of size $1$
in the definition of ER($1$) is not the reason why the moment assumption in Theorem~A is required.
\end{Remark}

\section{Results in  the noisy measurements setup}
\label{sec:noisy-measurements}
In previous sections, we considered the idealized scenario, in which the data was noiseless. Here, we will study the noisy setup: one observes $N$ couples $(z_i,X_i)_{i=1}^N$, and each $z_i$ is a noisy observation of $\inr{X_i,x_0}$:
\begin{equation}
  \label{eq:model}
z_i=\inr{X_i,x_0}+g_i, \quad i=1,\ldots,N.
\end{equation}
The goal is to obtain as much information as possible on the unknown vector $x_0$ with only the data $(z_i,X_i)_{i=1}^N$ at one's disposal, and for the sake of simplicity, we will assume that the $g_i$'s are independent Gaussian random variables $\cN(0,\sigma^2)$ that are also independent of the $X_i$'s.

Unlike the noiseless case, there is no hope of reconstructing $x_0$ from the given data, and instead of exact reconstruction, there are three natural questions that one may consider:
\begin{description}
\item{$\bullet$} The \textit{estimation problem}: given some norm $\norm{\cdot}$ on $\R^n$, one would like to construct a procedure $\hat x$ for which $\norm{\hat x-x_0}$ is as small as possible.
\item{$\bullet$} The \textit{prediction problem}: given a new (random, independent) `input' $X\in\R^n$, one has to find a good guess $\inr{\hat x,X}$ of the most likely associated output $z$, knowing that $(z,X)$ shares the same distribution with the other couples $(z_1,X_1),\ldots,(z_N,X_N)$.
\item{$\bullet$} The \textit{de-noising problem}: given a norm $\norm{\cdot}$ on $\R^N$ and a measurement matrix $\Gamma$, one has to construct $\hat x$ for which $\|\Gamma \hat x-\Gamma x_0\|$ is small.
\end{description}

These three problems are central in modern Statistics, and are featured in numerous statistical monographs, particularly in the context of the Gaussian regression model (Equation (\ref{eq:model})).

Recently, all three problems have been recast in a `high-dimensional'
scenario, in which the number of observations $N$ may be much smaller than the
ambient dimension $n$.  Unfortunately, such problems are often impossible to solve without additional assumptions, and just as in the noiseless case, the situation improves dramatically if $x_0$ has some \textit{low-dimensional structure}, for example, if it is $s$-sparse. The aim is therefore to design a procedure that performs as if the true dimension of the problem is $s$ rather than $n$, despite the noisy data.

To that end, $\ell_0$ penalization methods, sometimes called Model Selection procedures, have been introduced and studied extensively (see, e.g., \cite{MR2319879,MR1848946} for results in the context of the model (\ref{eq:model}), as well as in other examples). However, just as in the noise-free problem, the obvious downside of $\ell_0$ penalization methods is that they are not feasible computationally. This has lead to the introduction of convex relaxations, based on $\ell_1$ minimization.

Two well established $\ell_1$-based procedures are the \textit{LASSO} (see, e.g., \cite{MR1379242}) defined by
\begin{equation}
  \label{eq:lasso}
  \hat x_\lambda \in\argmin_{x\in\R^n}\Big(\frac{1}{N}\sum_{i=1}^N\big(z_i-\inr{X_i,x}\big)^2+\lambda \norm{x}_1\Big),
\end{equation}
and the Dantzig selector (see \cite{MR2382644}).

Both procedures may be implemented effectively, and their estimation and de-noising properties have been obtained under some assumptions on the measurement matrix (see, e.g. \cite{MR2807761,MR2533469,vdG07} or Chapters~7 and 8 in \cite{MR2829871}).

In this section, we shall focus on two such conditions on the measurement matrix. The first, called the Compatibility Condition, was introduced in \cite{vdG07} (see also Definition~2.1 in \cite{sara_decomposable}); the second, the Restricted Eigenvalue Condition, was introduced in \cite{MR2533469}.


\begin{Definition}\label{def:cc}
  Let $\Gamma$ be an ${N\times n}$ matrix. For $L>0$ and a set
  $S\subset\{1,\ldots,n\}$, the \textbf{compatibility constant}
  associated with $L$ and $S$ is
  \begin{equation}
    \label{eq:CC}
    \phi(L,S)=\sqrt{|S|}\min\Big(\norm{\Gamma \zeta_S-\Gamma\zeta_{S^c}}_{2}:\norm{\zeta_S}_1=1,\norm{\zeta_{S^c}}_1\leq L\Big),
  \end{equation}
  where $\zeta_S$ (resp. $\zeta_{S^c}$) denotes a vector that is supported in $S$ (resp. $S^c$).

$\Gamma$ satisfies the \textbf{Compatibility Condition for the set $S_0$} with constants $L>1$ and $c_0$ if $\phi(L,S_0)\geq c_0$; it satisfies the \textbf{uniform Compatibility Condition (CC) of order $s$} if $\min_{|S|\leq s}\phi(L,S)\geq c_0$.
\end{Definition}

A typical result for the LASSO in the Gaussian model~(\ref{eq:model})
and when $\Gamma$ satisfies the Compatibility Condition, is
Theorem~6.1 in \cite{MR2807761}:
\begin{Theorem}(\cite{MR2807761}, Theorem~6.1) Let $x_0\in\R^n$ be
  some fixed vector and assume that the data $(z_i,X_i)_{i=1}^N$ have
  been drawn according to the Gaussian regression
  model~\eqref{eq:model}. Denote by $\Gamma=N^{-1/2}\sum_{i=1}^N
  \inr{X_i,\cdot}f_i$ the measurement matrix.
Let $t>0$. If $S_0$ is the support of $x_0$ and $\lambda= 4 \sigma
\sqrt{(t^2+\log n)/N}$, then with probability larger than $1-2\exp(-t^2/2)$,
\begin{equation*}
  \norm{\Gamma \hat x_\lambda -\Gamma x_0 }_2^2 \leq
  \frac{64\sigma^2\norm{x_0}_0(t^2+ \log n)}{N \phi^2(3,S_0)} \mbox{ and } \norm{\hat
    x_\lambda - x_0}_1 \leq  \frac{64\sigma \norm{x_0}_0}{\phi^2(3,S_0)}\sqrt{\frac{t^2+\log n}{N}}.
\end{equation*}
\end{Theorem}

Even though the Compatibility Condition in $S_0$ suffices to show that the LASSO is an
effective procedure, the fact remains that $S_0$ is not known. And while a non-uniform approach is still possible (e.g., if $\Gamma$ is a random matrix, one may try showing that with high probability it satisfies the Compatibility Condition for the fixed, but unknown $S_0$), the uniform Compatibility Condition is a safer requirement -- and the one
we shall explore below.

Another uniform condition of a similar flavour is the Restricted Eigenvalue Condition from \cite{MR2533469}. To define it, let us introduce the following notation:  for $x\in\R^n$ and a set $S_0\subset\{1,\ldots,n\}$ of cardinality  $|S_0|\leq s$, let $S_1$ be the subset of indices of the $m$ largest coordinates of $(|x_i|)_{i=1}^n$ that are outside $S_0$. Let $x_{S_{01}}$ be the restriction of $x$ to the set $S_{01}=S_0\cup S_1$.

\begin{Definition}\label{def:rec}
  Let $\Gamma$ be an ${N\times n}$ matrix. Given $c_0\geq1$ and an integer $1\leq s\leq m\leq n$ for which $m+s\leq n$, the \textbf{restricted eigenvalue constant} is
  \begin{equation*}
\kappa(s,m,c_0) =    \min\Big(\frac{\norm{\Gamma x}_2}{\norm{x_{S_{01}}}_2}:S_0\subset\{1,\ldots,n\}, |S_0|\leq s, \norm{x_{S_0^c}}_1\leq c_0 \norm{x_{S_0}}_1\Big).
  \end{equation*}
The matrix $\Gamma$ satisfies the \textbf{Restricted Eigenvalue Condition (REC) of order $s$ with a constant $c$} if $\kappa(s,s,3)\geq c$.
\end{Definition}

Estimation and de-noising results follow from Theorem~6.1 (for the Dantzig selector) and Theorem~6.2 (for the LASSO) in \cite{MR2533469}, when the measurement matrix $\Gamma$, normalized by having the diagonal elements of $\Gamma^\top \Gamma$ equal $1$, satisfies the REC of an appropriate order and with a constant that is independent of the dimension. We also refer to Lemma~6.10 in \cite{MR2807761} for similar results that do not require normalization.

\vskip0.5cm

Because the two lead to bounds on the performance of the LASSO and the Dantzig selector, a question that comes to mind is whether there are matrices that satisfy the CC or the REC. And, as in Compressed Sensing, the only matrices that are known to satisfy those conditions for the optimal number of measurements (rows) are well-behaved random matrices (see \cite{MR2719855,MR3061256,Oliveira13,SaraMuro} for some examples).

Our aim in this final section is to extend our results to the noisy setup, by identifying almost necessary and sufficient moment assumptions for the CC and the REC. This turns out to be straightforward: on one hand, the proof of Theorem~A actually provides a stronger quantitative version of the exact reconstruction property; on the other, the uniform compatibility condition can be viewed as a \textit{quantitative version} of a geometric condition on the polytope $\Gamma B_1^n$ that characterizes Exact Reconstruction. A similar observation is true for the REC: it can be viewed as a quantitative version of the null space property (see \cite{MR1963681,MR1872845} and below) which is also equivalent to the exact reconstruction property.

\begin{Definition}
  Let $1\leq s \leq N$. A centrally symmetric polytope $P\subset\R^N$ is   \textbf{$s$-neighbourly} if every set of $s$ of its vertices, containing no antipodal pair, is the set of all vertices of some face of $P$.
\end{Definition}

It is well known \cite{donoho-polytope} that $\Gamma$ satisfies ER($s$) if and only if  $\Gamma B_1^n$ has $2n$ vertices and $\Gamma B_1^n$ is a centrally symmetric $s$-neighbourly polytope. It turns out that this property is characterized by the uniform CC.
\vskip0.3cm

\begin{Lemma}
Let $\Gamma$ be an ${N\times n}$ matrix. The following are equivalent:
  \begin{enumerate}
  \item $\Gamma B_1^n$ has $2n$ vertices and is $s$-neighbourly,
\item $\min\big(\phi(1,S):S\subset\{1,\ldots,n\}, |S|\leq s\big)>0$.
  \end{enumerate}
\end{Lemma}

In particular, $\min_{|S|\leq s}\phi(L,S)$ for some $L\geq1$ is a quantitative measure of the $s$-neighbourly property of $\Gamma B_1^n$: if $\Gamma B_1^n$ is $s$-neighbourly and has $2n$ vertices then the two sets
\begin{equation}\label{eq:sets_neighborly}
  \big\{\Gamma\zeta_S:\norm{\zeta_S}_1=1\big\} \mbox{ and } \big\{\Gamma \zeta_{S^c}:\norm{\zeta_{S^c}}_1\leq 1\big\}
\end{equation}
are  disjoint for every $|S|\leq s$. However, $\min_{|S|\leq
  s}\phi(1,S)$ measures how far the two sets are from one another,
uniformly over all subsets $S\subset\{1,\ldots,n\}$ of  cardinality at most $s$.
\vskip0.3cm
\proof Let $C_1,\ldots,C_n$ be the $n$ columns of $\Gamma$.
It follows from Proposition~2.2.13 and Proposition~2.2.16 in \cite{MR3113826} that $\Gamma B_1^n$ has $2n$ vertices and is a centrally symmetric $s$-neighbourly polytope if and only if for every $S\subset\{1,\ldots,n\}$ of cardinality $|S|\leq s$ and every choice of signs $(\eps_i)\in\{-1,1\}^S$,
\begin{equation}
  \label{eq:neighborly-1}
  {\rm conv}\big(\big\{\eps_i C_i:i\in S\big\}\big)\cap{\rm absconv}\big(\big\{C_j:j\notin S\big\}\big)=\emptyset.
\end{equation}
It is straightforward to verify that
\begin{equation*}
  \bigcup_{(\eps_i)\in\{\pm1\}^S} {\rm conv}\big(\big\{\eps_i C_i:i\in S\big\}\big)=\Big\{ \Gamma \zeta_S: \norm{\zeta_S}_1= 1\Big\}
\end{equation*}
and that
\begin{equation*}
   {\rm absconv}\big(\big\{C_j:j\notin S \big\}\big)=\Big\{ \Gamma \zeta_{S^c}: \norm{\zeta_{S^c}}_1\leq 1\Big\}.
\end{equation*}
As a consequence, (\ref{eq:neighborly-1}) holds for every $S\subset\{1,\ldots,n\}$ of cardinality at most $s$ if and only if  $\min\big(\phi(1,S):S\subset\{1,\ldots,n\}, |S|\leq s\big)>0$.
\endproof

An observation of a similar nature is true for the REC: it can be viewed as a quantitative measure of the null space property.

\begin{Definition}\label{def:nsp}
  Let $\Gamma$ be an ${N\times n}$ matrix. $\Gamma$ satisfies the \textbf{null space property of order $s$}  if it is invertible in the cone
  \begin{equation}
    \label{eq:cone_nsp}
    \big\{x\in\R^n: \mbox{ there exists } S\subset\{1,\ldots,n\}, |S|\leq s \mbox{ and } \norm{x_{S^c}}_1\leq \norm{x_{S}}_1\big\}.
  \end{equation}
\end{Definition}

In \cite{MR1963681,MR1872845}, the authors prove that $\Gamma$ satisfies ER($s$) if and only if it has the null space property of order $s$.

A natural way of quantifying the invertibility of $\Gamma$ in the cone \eqref{eq:cone_nsp} is to consider its smallest singular value, restricted to this cone, which is simply the REC $\kappa(s,n-s,1)$. Unfortunately, statistical properties of the LASSO and of the Dantzig selector are not known under the assumption that $\kappa(s, n-s, 1)$ is an absolute constant (though if $\kappa(s,s,3)$ is an absolute constant, LASSO is known to be optimal \cite{MR2533469}).

\vskip0.3cm

The main result of this section is the following:

\vskip0.3cm

\noindent {\bf Theorem E.} \textit{Let $L>0$, $1\leq s\leq n$ and $c_0>0$. Under the same assumptions as in Theorem~A and with the same probability estimate,  $\Gamma=N^{-1/2}\sum_{i=1}^N \inr{X_i,\cdot}f_i$ satisfies:
  \begin{enumerate}
  \item  A uniform compatibility condition of order $c_1 s $, namely that
  $$
      \min_{|S|\leq c_1s}\phi(L,S)\geq u^2\beta/4
  $$
  for $c_1=u^2\beta /(16e(1+L)^2)$.
\item A restricted eigenvalue condition of order $c_2 s$, with
$$
\kappa(c_2 s,m,c_0)\geq u^2\beta/4
$$
for any $1\leq m\leq n$, as long as $(1+c_0)^2 c_2\leq u^2 \beta /(16e)$.
  \end{enumerate}}

\textit{On the other hand, if $\Gamma$ is the matrix considered in Theorem~C, then with probability at least $1/2$, $\phi(1,\{e_1\})=0$ and $\kappa(1,m,1)=0$ for any $1\leq m\leq n$.}

\vskip0.3cm

Just like Theorem~A and Theorem~C, Theorem~E shows that the requirement that the coordinates of the measurement vector have $\log n$ moments is almost a necessary and sufficient condition for the uniform Compatibility Condition and the Restricted Eigenvalue Condition to hold. Moreover, it shows the significance of the small-ball condition, even in the noisy setup.

It also follows from Theorem~E that if $X$ satisfies the small-ball condition and its coordinates have $\log n$ well-behaved moments as in Theorem~A, then $\Gamma B_1^n$ has $2n$ vertices and is $s$-neighbourly with high probability for $N\sim s \log(en/s)$. In particular, this improves Theorem~4.3 in \cite{MR2796091} by a logarithmic factor for matrices generated by subexponential variables.

\textbf{Proof of Theorem~E:} Fix a constant $c_1$ to be named later and let  $S\subset\{1,\ldots,n\}$ of cardinality $|S|\leq c_1 s$. Let $\zeta_S\in\R^n$ be a vector supported on $S$ with $\norm{\zeta_S}_1=1$ and let $\zeta_{S^c}\in\R^n$ be supported on $S^c$ with $\norm{\zeta_{S^c}}_1\leq L$.

Consider $\gamma=(\zeta_S-\zeta_{S^c})/\norm{\zeta_S-\zeta_{S^c}}_2$. Since
$$
\norm{\zeta_S-\zeta_{S^c}}_2 \geq \norm{\zeta_S}_2 \geq \frac{\norm{\zeta_S}_1}{\sqrt{|S|}}= \frac{1}{\sqrt{|S|}},
$$
it follows that
$\gamma\in \big((1+L)\sqrt{|S|}\big)B_1^n \cap S^{n-1}$.

 Recall that by (\ref{eq:in-proof-cone}), if $r=(1+L)^2 c_1 s \leq s u^2\beta/(16e)$, then $\norm{\Gamma \gamma}_2\geq (u^2 \beta)/4$. Therefore,
\begin{equation*}
  \norm{\Gamma \zeta_S-\Gamma\zeta_{S^c}}_2\geq \frac{u^2\beta}{4}\norm{\zeta_S-\zeta_{S^c}}_2\geq \frac{u^2\beta}{4}\norm{\zeta_S}_2\geq \frac{u^2 \beta\norm{\zeta_S}_1}{4\sqrt{|S|}}=\frac{u^2 \beta}{4\sqrt{|S|}},
\end{equation*}
and thus $\min_{|S|\leq c_1s}\phi(L,S)\geq u^2\beta/4$ for $c_1 =u^2\beta/\big(16 e (1+L)^2 \big)$.

Turning to the REC, fix a constant $c_2$ to be named later. Consider $x$ in the cone and let $S_0\subset\{1,\ldots,n\}$ of cardinality $|S_0|\leq c_2s$ for which $\norm{x_{S_0^c}}_1\leq c_0 \norm{x_{S_0}}_1$. Let $S_1\subset\{1,\ldots,n\}$ be the set of indices of the $m$ largest coordinates of $(|x_i|)_{i=1}^n$ that are outside  $S_0$ and put $S_{01}=S_0\cup S_1$.

Observe that $\norm{x}_1  \leq (1+c_0)\norm{x_{S_0}}_1 \leq (1+c_0)\sqrt{|S_0|}\norm{x}_2$; hence $x/\norm{x}_2\in \big((1+c_0)\sqrt{|S_0|}\big)B_1^n\cap S^{n-1}$. Applying (\ref{eq:in-proof-cone}) again, if $(1+c_0)^2c_2 s \leq s u^2 \beta/(16 e)$, then $\norm{\Gamma x}_2\geq \big((u^2\beta)/4\big)\norm{x}_2$. Thus,
\begin{equation*}
  \frac{\norm{\Gamma x}_2}{\norm{x_{S_{01}}}_2}\geq \frac{\norm{\Gamma x}_2}{\norm{x}_2}\geq \frac{u^2 \beta}{4}
\end{equation*}
and $\kappa(c_2 s,m,c_0)\geq u^2\beta/4$ for any $1\leq m\leq n$, as long as $(1+c_0)^2 c_2\leq u^2 \beta /(16e)$.

The proof of the second part of Theorem~E is an immediate corollary of the construction used in Theorem~C. Recall that with probability at least $1/2$, $\Gamma e_1 \in {\rm absconv}(\Gamma e_j : j \in \{2,...,n\})$. Setting $J=\{e_2,...,e_n\}$, there is $\zeta \in B_1^J$ for which $\norm{\Gamma e_1-\Gamma \zeta}_2=0$. Therefore, $\phi(1,\{e_1\})=0$ and $\kappa(1,m,1)=0$ for any $1\leq m\leq n$, as claimed.
\endproof

\begin{Remark} The results obtained in Theorem~A and in parts (1) and (2) of Theorem~E are also valid for the
  normalized (columns wise) measurement matrix:
\begin{equation*}
  {\Gamma}_1 =\Gamma \tilde D^{-1} \mbox{ where } \tilde D = {\rm diag}\big(\norm{\Gamma e_1}_2,\ldots,\norm{\Gamma e_n}_2\big).
\end{equation*}
The proof is almost identical to the one used for $\Gamma$ itself, even though  $\Gamma_1$ does not have independent rows vectors, due to the normalization. For the sake of brevity, we will not present the straightforward proof of this observation.

Finally, the counterexample constructed in the proof of Theorem~C and in which a typical $\Gamma$ does not satisfy ER($1$), does not necessarily generate $\Gamma B_1^n$ that is not $s$-neighbourly. Indeed, an inspection of the construction shows that the reason ER($1$) fails is that $\Gamma B_1^n$ has less than $2n-2$ vertices, rather than that $\Gamma B_1^n$ is not $s$-neighbourly. Thus, the question of whether a moment condition is necessary for the random polytope $\Gamma B_1^n$ to be $s$-neighbourly with probability at least $1/2$ is still unresolved.

\end{Remark}

\begin{footnotesize}

\bibliographystyle{plain}

\bibliography{biblio}

\def\cprime{$'$}
\begin{thebibliography}{10}

\bibitem{MR2796091}
Rados{\l}aw Adamczak, Alexander~E. Litvak, Alain Pajor, and Nicole
  Tomczak-Jaegermann.
\newblock Restricted isometry property of matrices with independent columns and
  neighborly polytopes by random sampling.
\newblock {\em Constr. Approx.}, 34(1):61--88, 2011.

\bibitem{MR840631}
Keith Ball.
\newblock Cube slicing in {${\bf R}^n$}.
\newblock {\em Proc. Amer. Math. Soc.}, 97(3):465--473, 1986.

\bibitem{MR2533469}
Peter~J. Bickel, Ya'acov Ritov, and Alexandre~B. Tsybakov.
\newblock Simultaneous analysis of {L}{A}{S}{S}{O} and {D}antzig selector.
\newblock {\em Ann. Statist.}, 37(4):1705--1732, 2009.

\bibitem{MR1848946}
Lucien Birg{\'e} and Pascal Massart.
\newblock Gaussian model selection.
\newblock {\em J. Eur. Math. Soc. (JEMS)}, 3(3):203--268, 2001.

\bibitem{BouLugMass13}
St{\'e}phane Boucheron, G{\'a}bor Lugosi, and Pascal Massart.
\newblock {\em Concentration Inequalities: A Nonasymptotic Theory of
  Independence.}
\newblock Oxford University Press, 2013.
\newblock ISBN 978-0-19-953525-5.

\bibitem{MR2807761}
Peter B{\"u}hlmann and Sara~A. van~de Geer.
\newblock {\em Statistics for high-dimensional data}.
\newblock Springer Series in Statistics. Springer, Heidelberg, 2011.
\newblock Methods, theory and applications.

\bibitem{C:08}
Emmanuel~J. Cand{\`e}s.
\newblock The restricted isometry property and its implications for compressed
  sensing.
\newblock {\em C. R. Math. Acad. Sci. Paris}, 346(9-10):589--592, 2008.

\bibitem{MR2236170}
Emmanuel~J. Cand{\`e}s, Justin Romberg, and Terence Tao.
\newblock Robust uncertainty principles: exact signal reconstruction from
  highly incomplete frequency information.
\newblock {\em IEEE Trans. Inform. Theory}, 52(2):489--509, 2006.

\bibitem{MR2230846}
Emmanuel~J. Cand{\`e}s, Justin~K. Romberg, and Terence Tao.
\newblock Stable signal recovery from incomplete and inaccurate measurements.
\newblock {\em Comm. Pure Appl. Math.}, 59(8):1207--1223, 2006.

\bibitem{MR2300700}
Emmanuel~J. Cand{\`e}s and Terence Tao.
\newblock Near-optimal signal recovery from random projections: universal
  encoding strategies?
\newblock {\em IEEE Trans. Inform. Theory}, 52(12):5406--5425, 2006.

\bibitem{MR2382644}
Emmanuel~J. Cand{\`e}s and Terence Tao.
\newblock The {D}antzig selector: statistical estimation when {$p$} is much
  larger than {$n$}.
\newblock {\em Ann. Statist.}, 35(6):2313--2351, 2007.

\bibitem{MR3113826}
Djalil Chafa{\"{\i}}, Olivier Gu{\'e}don, Guillaume Lecu{\'e}, and Alain Pajor.
\newblock {\em Interactions between compressed sensing random matrices and high
  dimensional geometry}, volume~37 of {\em Panoramas et Synth\`eses [Panoramas
  and Syntheses]}.
\newblock Soci\'et\'e Math\'ematique de France, Paris, 2012.

\bibitem{MR1854649}
Scott~Shaobing Chen, David~L. Donoho, and Michael~A. Saunders.
\newblock Atomic decomposition by basis pursuit.
\newblock {\em SIAM Rev.}, 43(1):129--159, 2001.
\newblock Reprinted from SIAM J. Sci. Comput. {{\bf{2}}0} (1998), no. 1, 33--61
  (electronic) [ MR1639094 (99h:94013)].

\bibitem{claerbout}
Jon~F. Claerbout and Francis Muir.
\newblock Robust modeling with erratic data.
\newblock {\em Geophysics}, 38(5):826--844, 1973.

\bibitem{MR1666908}
V{\'{\i}}ctor~H. de~la Pe{\~n}a and Evarist Gin{\'e}.
\newblock {\em Decoupling}.
\newblock Probability and its Applications (New York). Springer-Verlag, New
  York, 1999.
\newblock From dependence to independence, Randomly stopped processes.
  $U$-statistics and processes. Martingales and beyond.

\bibitem{donoho-polytope}
David~L. Donoho.
\newblock Neighborly polytopes and sparse solutions of under-determined linear
  equations.
\newblock Technical report, Department of Statistics, Standford University,
  2005.

\bibitem{MR2241189}
David~L. Donoho.
\newblock Compressed sensing.
\newblock {\em IEEE Trans. Inform. Theory}, 52(4):1289--1306, 2006.

\bibitem{MR1963681}
David~L. Donoho and Michael Elad.
\newblock Optimally sparse representation in general (nonorthogonal)
  dictionaries via {$l^1$} minimization.
\newblock {\em Proc. Natl. Acad. Sci. USA}, 100(5):2197--2202 (electronic),
  2003.

\bibitem{MR1872845}
David~L. Donoho and Xiaoming Huo.
\newblock Uncertainty principles and ideal atomic decomposition.
\newblock {\em IEEE Trans. Inform. Theory}, 47(7):2845--2862, 2001.

\bibitem{MR1154788}
David~L. Donoho and Ben~F. Logan.
\newblock Signal recovery and the large sieve.
\newblock {\em SIAM J. Appl. Math.}, 52(2):577--591, 1992.

\bibitem{MR512411}
R.~M. Dudley.
\newblock Central limit theorems for empirical measures.
\newblock {\em Ann. Probab.}, 6(6):899--929 (1979), 1978.

\bibitem{MR3134335}
Simon Foucart.
\newblock Stability and robustness of {$\ell_1$}-minimizations with {W}eibull
  matrices and redundant dictionaries.
\newblock {\em Linear Algebra Appl.}, 441:4--21, 2014.

\bibitem{MR3100033}
Simon Foucart and Holger Rauhut.
\newblock {\em A mathematical introduction to compressive sensing}.
\newblock Applied and Numerical Harmonic Analysis. Birkh\"auser/Springer, New
  York, 2013.

\bibitem{MR2829871}
Vladimir Koltchinskii.
\newblock {\em Oracle inequalities in empirical risk minimization and sparse
  recovery problems}, volume 2033 of {\em Lecture Notes in Mathematics}.
\newblock Springer, Heidelberg, 2011.
\newblock Lectures from the 38th Probability Summer School held in Saint-Flour,
  2008, {\'E}cole d'{\'E}t{\'e} de Probabilit{\'e}s de Saint-Flour.
  [Saint-Flour Probability Summer School].

\bibitem{Shahar-Vladimir}
Vladimir Koltchinskii and Shahar Mendelson.
\newblock Bounding the smallest singular value of a random matrix without
  concentration.
\newblock Technical report, Technion and Georgia Tech, 2013.
\newblock arXiv:1312.3580.

\bibitem{MR1457628}
Rafa{\l} Lata{\l}a.
\newblock Estimation of moments of sums of independent real random variables.
\newblock {\em Ann. Probab.}, 25(3):1502--1513, 1997.

\bibitem{LT:91}
Michel Ledoux and Michel Talagrand.
\newblock {\em Probability in {B}anach spaces}, volume~23 of {\em Ergebnisse
  der Mathematik und ihrer Grenzgebiete (3) [Results in Mathematics and Related
  Areas (3)]}.
\newblock Springer-Verlag, Berlin, 1991.
\newblock Isoperimetry and processes.

\bibitem{logan}
Ben Logan.
\newblock {\em Properties of High-Pass Signals}.
\newblock PhD thesis, Columbia University, New York, 1965.

\bibitem{MR2319879}
Pascal Massart.
\newblock {\em Concentration inequalities and model selection}, volume 1896 of
  {\em Lecture Notes in Mathematics}.
\newblock Springer, Berlin, 2007.
\newblock Lectures from the 33rd Summer School on Probability Theory held in
  Saint-Flour, July 6--23, 2003, With a foreword by Jean Picard.

\bibitem{Shahar-COLT}
Shahar Mendelson.
\newblock Learning without concentration.
\newblock Journal of the ACM. To appear.

\bibitem{Shahar-Gelfand}
Shahar Mendelson.
\newblock A remark on the diameter of random sections of convex bodies.
\newblock Geometric Aspects of Functional Analysis (GAFA Seminar Notes).
  Lecture notes in Mathematics 2116, pages 395--404, 2014.

\bibitem{MR2453368}
Shahar Mendelson, Alain Pajor, and Nicole Tomczak-Jaegermann.
\newblock Uniform uncertainty principle for {B}ernoulli and subgaussian
  ensembles.
\newblock {\em Constr. Approx.}, 28(3):277--289, 2008.

\bibitem{MR1320206}
B.~K. Natarajan.
\newblock Sparse approximate solutions to linear systems.
\newblock {\em SIAM J. Comput.}, 24(2):227--234, 1995.

\bibitem{Oliveira13}
Roberto~Imbuzeiro Oliveira.
\newblock The lower tail of random quadratic forms, with applications to
  ordinary least squares and restricted eigenvalue properties.
\newblock Technical report, IMPA, 2013.
\newblock arXiv:1312.2903.

\bibitem{MR2719855}
Garvesh Raskutti, Martin~J. Wainwright, and Bin Yu.
\newblock Restricted eigenvalue properties for correlated {G}aussian designs.
\newblock {\em J. Mach. Learn. Res.}, 11:2241--2259, 2010.

\bibitem{MR890930}
B.~A. Rogozin.
\newblock An estimate for the maximum of the convolution of bounded densities.
\newblock {\em Teor. Veroyatnost. i Primenen.}, 32(1):53--61, 1987.

\bibitem{MR2417886}
Mark Rudelson and Roman Vershynin.
\newblock On sparse reconstruction from {F}ourier and {G}aussian measurements.
\newblock {\em Comm. Pure Appl. Math.}, 61(8):1025--1045, 2008.

\bibitem{RV13}
Mark Rudelson and Roman Vershynin.
\newblock Small ball probabilities for linear images of high dimensional
  distributions.
\newblock Technical report, University of Michigan, 2013.

\bibitem{MR3061256}
Mark Rudelson and Shuheng Zhou.
\newblock Reconstruction from anisotropic random measurements.
\newblock {\em IEEE Trans. Inform. Theory}, 59(6):3434--3447, 2013.

\bibitem{MR857796}
Fadil Santosa and William~W. Symes.
\newblock Linear inversion of band-limited reflection seismograms.
\newblock {\em SIAM J. Sci. Statist. Comput.}, 7(4):1307--1330, 1986.

\bibitem{MR3127875}
Nikhil Srivastava and Roman Vershynin.
\newblock Covariance estimation for distributions with {$2+\varepsilon$}
  moments.
\newblock {\em Ann. Probab.}, 41(5):3081--3111, 2013.

\bibitem{tbm}
Howard~L. Taylor, Stephen~C. Banks, and John~F. McCoy.
\newblock Deconvolution with the l1-norm.
\newblock {\em Geophysics}, 44(1):39--52, 1979.

\bibitem{MR1379242}
Robert Tibshirani.
\newblock Regression shrinkage and selection via the lasso.
\newblock {\em J. Roy. Statist. Soc. Ser. B}, 58(1):267--288, 1996.

\bibitem{vdG07}
Sara~A. van~de Geer.
\newblock The deterministic lasso.
\newblock {\em In JSM proceedings. American Statistical Association}, 140,
  2007.

\bibitem{sara_decomposable}
Sara~A. van~de Geer.
\newblock Weakly decomposable regularization penalties and structured sparsity.
\newblock Technical report, ETH Z{\"u}rich, 2013.

\bibitem{SaraMuro}
Sara~A. van~de Geer and Alan Muro.
\newblock On higher order isotropy conditions and lower bounds for sparse
  quadratic forms.
\newblock Technical report, ETH Z{\"u}rich, 2014.

\bibitem{vanderVaartWellner}
Aad~W. van~der Vaart and Jon~A. Wellner.
\newblock {\em Weak convergence and empirical processes}.
\newblock Springer Series in Statistics. Springer-Verlag, New York, 1996.
\newblock With applications to statistics.

\end{thebibliography}
\end{footnotesize}

\end{document}